\documentclass[11pt]{amsart}
\usepackage{latexsym,amssymb,amsmath,hyperref}
\usepackage{amsthm, dsfont}
\usepackage[all]{xy}
\usepackage[short,nodayofweek]{datetime}
\usepackage[active]{srcltx}

\def\GG{{\mathbb G}}
\def\NN{{\mathbb N}}
\def\ZZ{{\mathbb Z}}
\def\QQ{{\mathbb Q}}

\def\WW{{\mathbb W}}
\def\KK{{\mathbb K}}

\def\O{{\mathcal O}}

\newcommand{\kk}{{\Bbbk}}

\newcommand{\Ga}{{\GG_a}}

\newcommand{\thetaj}{\theta^{(j)}}
\newcommand{\thetak}{\theta^{(k)}}

\newcommand{\thetai}{\theta^{(i)}}
\newcommand{\thetan}{\theta^{(n)}}

\def\id{{\rm id}}

\def\isom{\cong}
\def\congr{\equiv}

\def\ideal{\unlhd}

\DeclareMathOperator{\Mat}{Mat}

\theoremstyle{plain}
\newtheorem{thm}{Theorem}[section]

\newtheorem{lem}[thm]{Lemma}
\newtheorem{prop}[thm]{Proposition}

\theoremstyle{definition}

\theoremstyle{remark}
  \newtheorem{rem}[thm]{Remark}

\newtheoremstyle{Acknowledgements}
  {}
    {}
     {}
     {}
    {\bfseries}
    {}
     {.5em}
     {\thmname{#1}\thmnumber{ }\thmnote{ (#3)}}
\theoremstyle{Acknowledgements}
\newtheorem{ack}{Acknowledgements.}

\date{\today} 

\begin{document}

\title[Finite generation of $\Ga$-invariants]{On the finite generation of additive group invariants in positive characteristic}

\author{Emilie Dufresne}
\address{\rm {\bf Emilie Dufresne}, Mathematics Center of Heidelberg
  (MATCH), Heidelberg University, Im Neuenheimer Feld 368, 69120 Heidelberg, Germany }
\curraddr{}
\email{\sf emilie.dufresne@iwr.uni-heidelberg.de}
\thanks{}

\author{Andreas Maurischat}
\address{\rm {\bf Andreas Maurischat}, Interdisciplinary Center for
  Scientific Computing (IWR), Heidelberg 
University, Im Neuenheimer Feld 368, 69120 Heidelberg, Germany }
\email{\sf andreas.maurischat@iwr.uni-heidelberg.de}



\begin{abstract}
 Roberts, Freudenburg, and Daigle and Freudenburg have given the smallest
counterexamples to Hilbert's fourteenth problem as rings of invariants of
algebraic groups. Each is of an action of the additive group on a
finite dimensional vector space over a field of characteristic zero, and thus,
each is the kernel of a locally nilpotent derivation. In positive
characteristic, additive group actions correspond to locally finite iterative
higher derivations. We set up characteristic-free analogs of the three examples,
and show that, contrary to characteristic zero, in every positive charateristic,
the invariants are finitely generated.
\end{abstract}

\maketitle



\section{Introduction}

A main topic of interest in Invariant Theory is the question of the finite
generation of invariant rings. Namely, if $B$ is a finitely generated algebra
over a field $\kk$, and $G$ is a group acting on $B$ via $\kk$-algebra
homomorphisms, one asks if the ring of invariants $B^G$ is finitely
generated. This is a special case of Hilbert's Fourteenth Problem. When $G$ is a
finite group, this question has a positive answer for any $B$ (from the work of
Hilbert and Noether \cite{em-deig,em-deigcp}). If $G$ is an affine algebraic group
 acting regularly on $B$, then whenever $G$ is  reductive, the ring of invariants is finitely generated (cf. Nagata \cite{mn:igar}). On the other hand, if
$G$ is not reductive, then Popov (cf. \cite{vlp:ohti}) showed that there exists
a finitely generated $\kk$-algebra $B$  and an action of $G$ on $B$ such that
$B^G$ is not finitely generated. These results are valid in arbitrary
characteristic.

Consider actions on a polynomial ring $B=\kk[x_1,\dots, x_n]$ of the simplest
non reductive group, the additive group $\Ga$. In characteristic zero, several
facts are known. By the Maurer-Weitzen\-b\"ock Theorem (cf. \cite{rw:uilg}),
every linear action of $\Ga$ on $B$ has finitely generated
invariants. If $n\leq 3$, then $B^\Ga$ is finitely generated for any algebraic
action (cf. \cite{oz:iaqph}). In dimensions $5$, $6$ and $7$, however, there are
well-known counterexamples. For $n=7$, Roberts (cf.
\cite{pr:aigsbpsrnchfp,jkd-drf:gac3c7}) gave an example of a $\Ga$-action where
$B^\Ga$ is not finitely generated (from now on refered to as (R7)). For $n=6$
and $5$, examples where constructed by Freudenburg (F6) (cf. \cite{gf:achfpd6}), and Daigle and
Freudenburg (DF5) (cf. \cite{dd-gf:achfpd5}). These three counterexamples can be used to construct
counterexamples in any dimension $n\geq 5$
(cf. \cite{ae-tj:ked}). The dimension $4$ case is still open for general
$\Ga$-actions. For more information on $\Ga$-actions in characteristic
zero, we refer the reader to  the excellent book of Freudenburg \cite{gf:atlnd}.

In positive characteristic, much less is known. It is known that $B^\Ga$ is
finitely generated  if $n\leq 3$ (cf. \cite{oz:iaqph}), and even polynomial if
$n\leq 2$ (cf. \cite{mm:gaaap}). The finite generation of the invariants was
also proved for special classes of linear actions of $\Ga$ (cf. e.g.
\cite{css:otwit}, \cite{rt:opirc1Gam}, \cite{rt:aacklfihd}), but not for all
linear actions. On the other
hand, to the authors's knowledge, there is no example of an algebraic
$\Ga$-action on a polynomial ring where the ring of invariants was proven not to
be finitely generated.

Locally finite iterative higher derivations (lfihd) are a generalization of
locally nilpotent derivations which behave well in all characteristics: an
algebraic action of the additive group always corresponds to a lfihd. In this
paper, we adopt this point of view, and consider the positive characteristic
analogs of the counterexamples (R7), (F6), and (DF5) mentioned above. Our main
result is that, contrary to characteristic zero, in every positive
characteristic, the arising rings of invariants are finitely generated. For
Robert's example (R7), this was done by Kurano (c.f. \cite{kk:pcfgsrarcfph}),
but not for (F6) and (DF5). Furthermore, our approach is both constructive and
more elementary.

The paper is structured as follows. In Section \ref{sec:setup}, we recall some
facts concerning lfihd, and set up the examples. In Section~\ref{sec:main-results}, we prove that the invariants are finitely generated.  We
use the fact that the examples are related by homomorphisms which respect the
lfihd. The key of the argument is the existence of a special invariant. We
deduce the existence of the special element for (F6) and (R7) from the
existence of the special element for (DF5). In Section~\ref{sec:dim5},
we establish the existence of the special element for (DF5). Our rather
technical construction is in the spirit of van den Essen's simpler proof for
(DF5) in characteristic zero (cf. \cite{ae:asshfp}). Finally in Section~\ref{sec:generalization}, we explain how our argument can be extended to more general versions of (DF5), (F6), and (R7).

\begin{ack}
 We thank Florian Heiderich for giving us the idea to explore the correspondence
between additive group actions and locally finite iterative higher derivations.
We also thank Kazuhiko Kurano for calling our attention to his paper, and thus
encouraging us to extend our argument to more general versions of (F6) and
(DF5).
\end{ack}



\section{Setup of the Examples}\label{sec:setup}

We first review a purely algebraic description of $\Ga$-actions.
Throughout this paper, $\kk$ denotes an
algebraically closed field.\footnote{This is not a restriction: the property of
finite generation of invariants is stable under algebraic extensions of the base
field.} Let $B$ be a finitely generated $\kk$-algebra. An algebraic action
of $\Ga$ on $B$ is uniquely determined by a $\kk$-algebra homomorphism
$\theta:B\to B\otimes_{\kk} \kk[U]=B[U]$, where $\kk[U]\isom
\kk[\Ga]$ is the ring of regular functions on $\Ga$. The correspondence is given by $\sigma\cdot b=\theta(b)|_{U=\sigma}$ for all $b\in B$ and $\sigma\in \Ga(\kk)=\kk$. Any $\kk$-algebra homomorphism $\theta:B\to B[U]$ defines a family of $\kk$-linear maps $\left(\thetan\right)_{n\geq 0}$ via
$$\theta(b)=:\sum_{n=0}^\infty \thetan(b) U^n$$
for all $b\in B$.  A family $\left(\thetan\right)_{n\geq 0}$ (resp. $\theta$) corresponds to a
$\Ga$-action if and only if it fulfills the following properties (cf. \cite{mm:ariihd}):
\begin{enumerate}
\item[(1)] $\theta^{(0)}=\id_B,$
\item[(2)] for all $n\geq 0$ and $a,b\in B$, one has $\thetan(ab)=\sum_{i+j=n}
\thetai(a) \thetaj(b)$,
\item[(3)] for all $b\in B$, there is $n\geq 0$ such that $\thetaj(b)=0$ for
all $j\geq n$,
\item [(4)] for all $j,k\geq 0$ and $b\in B$, one has
$\thetaj(\thetak(b))=\binom{j+k}{j}\theta^{(j+k)}(b)$.
\end{enumerate}

Whereas Properties (2) and (3) are equivalent to $\theta$ being a $\kk$-algebra homomorphism, Properties (1) and (4) ensure that $\theta$ really determines a $\Ga$-action. A
family $\left(\thetan\right)_{n\geq 0}$ fulfilling these properties (resp. the
corresponding $\theta$) is called a {\em locally finite iterative higher
derivation} (lfihd) on $B$. Throughout this paper, we adopt the point of view of lfihd. 

Note that $\theta^{(1)}$ is a derivation in the usual sense, and in
characteristic zero, Property (4) implies that $\thetan=\frac{1}{n!}
(\theta^{(1)})^n$. Thus, in characteristic zero, lfihd are in one to one
correspondence with those derivations which are locally nilpotent (by Property (3)).
Accordingly, the $\kk$-algebra $B^\theta:=\{ b\in B \mid \thetan(b)=0\; \forall n\geq 1\}$ is
often called the ring of constants of $\theta$. Since, it coincides with the ring of
invariants $B^\Ga$ of the corresponding $\Ga$-action, we refer to
$B^\theta$ as the ``ring of invariants''.

In characteristic zero,  (DF5), (F6) and (R7) are realized as $\kk$-algebras with a locally nilpotent derivation:
the counterexample of Daigle and Freudenburg (DF5) is given by
$$\kk[x,s,t,u,v]\quad \text{and} \quad\theta^{(1)}=x^3
\frac{\partial}{\partial s} + s \frac{\partial}{\partial t}+t
\frac{\partial}{\partial u}+x^2 \frac{\partial}{\partial v},$$
the counterexample of Freudenburg  (F6) is given by
$$\kk[x,y,s,t,u,v]\quad \text{and} \quad\theta^{(1)}=x^3
\frac{\partial}{\partial s} + y^3 s \frac{\partial}{\partial t}+y^3 t
\frac{\partial}{\partial u}+x^2 y^2 \frac{\partial}{\partial v},$$
and the counterexample of Roberts (R7) is given by
$$\kk[x_1,x_2,x_3,y_1,y_2,y_3,v]\ \text{and}
\ \theta^{(1)}=x_1^3 \frac{\partial}{\partial y_1} + x_2^3 
\frac{\partial}{\partial y_2}+x_3^3  \frac{\partial}{\partial y_3}+x_1^2
x_2^2 x_3^2 \frac{\partial}{\partial v}.$$

To obtain locally finite iterative higher derivations which make sense in all
positive characteristics, we rescale the variables $t$ and $u$  in (DF5) and
(F6) by a factor of $2$ and $6$, respectively. Characteristic-free formulations
of the examples are therefore given by:

\noindent {\bf Daigle-Freudenburg's example} (DF5): \\$B_5:=\kk[x,s,t,u,v]$,
\begin{align*}
\theta(x)&=x,& \theta(s)&=s+x^3 U, \\
\theta(t)&=t+2sU+x^3U^2,& \theta(u)&=u+3tU+3sU^2+x^3U^3,\\ 
\theta(v)&=v+x^2 U. && 
\end{align*} 
{\bf Freudenburg's example} (F6): \\$B_6:=\kk[x,y,s,t,u,v]$,
\begin{align*} 
\theta(x)&=x,& \theta(y)&=y,\\
\theta(s)&=s+x^3 U,& \theta(t)&=t+2y^3sU+x^3y^3U^2,\\ 
\theta(u)&=u+3y^3tU+3y^6sU^2+x^3y^6U^3,& \theta(v)&=v+x^2y^2 U.
\end{align*}
{\bf Roberts's example} (R7):\\ $B_7:=\kk[x_1,x_2,x_3,y_1,y_2,y_3,v]$,
\begin{align*}
\theta(x_i)&=x_i,& \theta(y_i)&=y_i+x_i^3 U\quad (i=1,2,3),\\ 
\theta(v)&=v+ x_1^2 x_2^2 x_3^2 U.&&
\end{align*}

On the $\kk$-algebras $B_5, B_6, B_7$, we define gradings $w_5,w_6,w_7$ by assigning the following degrees:
\[\begin{array}{l}
w_5(x)=1,~w_5(s)=w_5(t)=w_5(u)=3,~ w_5(v)=2;\\
w_6(x)=w_6(y)=1,~w_6(s)=3,~w_6(t)=6,~w_6(u)=9, w_6(v)=4;\\
w_7(x_i)=1,~w_7(y_i)=3,~(i=1,2,3),~w_7(v)=6.\\
  \end{array}\]
With respect to these gradings, the lfihd and the corresponding $\Ga$-actions
are homogeneous, and so the rings of invariants are graded subalgebras. This
provides useful additional structure. We will also use an additional grading
$w_4$ on $B_5$ which is given by:
$$w_4(x)=0,~w_4(s)=1,~w_4(t)=2,~w_4(u)=3,~w_4(v)=1.$$

We now have the proper setup to state our main theorem:

\begin{thm}\label{thm:MainThm}
In every positive characteristic, the rings of invariants $B_5^\theta$,
$B_6^\theta$, and $B_7^\theta$, are finitely generated.
\end{thm}



\section{Main Results}\label{sec:main-results}

This section presents the main steps of our argument to prove
Theorem \ref{thm:MainThm}. We will make use of Theorem
\ref{thm:special-invariant-B5}, which states the existence of a certain
invariant in $B_5$.
First, we describe the connection between the examples:

\begin{lem}\label{lem:homomorphisms}~
\begin{enumerate}
\item The ring $B_5$ is isomorphic to $B_6/(y-1)$ and the lfihd on $B_5$ is the
lfihd induced by this isomorphism from the lfihd on $B_6$.
\item A homomorphism $\alpha:B_6\to B_7$ respecting the lfihd is
given by:
\begin{align*}
\alpha(x)&=x_1,& \alpha(y)&=x_2 x_3,\\
\alpha(s)&=y_1,&¸\alpha(t)&=(x_3^3y_2+x_2^3y_3)y_1-x_1^3y_2y_3,\\
\alpha(v)&=v, &
\alpha(u)&=(x_2^6y_3^2+x_2^3x_3^3y_2y_3+x_3^6y_2^2)y_1-
(x_3^3y_2+x_2^3y_3)x_1^3y_2y_3
\end{align*}
\end{enumerate}
\begin{proof}
This can be verified by a short computation.
\end{proof}
\end{lem}

\begin{prop}\label{prop:special-invariant}
Let $\kk$ be of positive characteristic $p$. Then in each of $B_5$, $B_6$, and
$B_7$, there exists a homogeneous invariant of the form $v^p+vb'-b$ such
that $v$ does not appear in $b$ and $b'$.
\begin{proof}
By Theorem \ref{thm:special-invariant-B5}, there is a $w_5$-homogeneous element
$v^p+vb'-b\in B_5^\theta$, such that $w_4(b)=p$ and $w_4(b')=p-1$ . Since $w_5(v^p)=2p$, we have
$w_5(b)=2p$ and $w_5(b')=2p-2$. Furthermore, $B_5\isom B_6/(y-1)$, and so we obtain a 
$w_6$-homogeneous invariant element $v^p+v\tilde{b'}-\tilde{b}$ in
$B_6[\frac{1}{y}]$ by
homogenizing $b$ and $b'$ inside $B_6[\frac{1}{y}]$, that is, by setting 
$$\tilde{b}:=y^{4p}\cdot
b\left(\frac{x}{y},\frac{s}{y^3},\frac{t}{y^6},\frac{u}{y^9}\right)\textrm{  and
~~ }\tilde{b'}:=y^{4p-4}\cdot
b'\left(\frac{x}{y},\frac{s}{y^3},\frac{t}{y^6},\frac{u}{y^9}\right).$$
Denote by $d_x,d_s,d_t,d_u$ the exponents in a
monomial of $b$ of the variables $x,s,t,u$, respectively. By the conditions on
$b$, we have
$2p=w_5(b)=d_x+3d_s+3d_t+3d_u$, and $p=w_4(b)=d_s+2d_t+3d_u$.
It follows that
$$d_x+3d_s+6d_t+9d_u\leq d_x+3d_s+3d_t+3d_u+2d_s+4d_t+6d_u=2p+2p=4p.$$
Hence, the exponent of $y$ in the denominator of
$b(\frac{x}{y},\frac{s}{y^3},\frac{t}{y^6},\frac{u}{y^9})$ is less or equal to
$4p$, and so $\tilde{b}\in B_6$. A similar argument shows that $\tilde{b'}\in B_6$.
Therefore, $v^p+v\tilde{b'}-\tilde{b}$ is an invariant element in $B_6$.

Finally, applying the homomorphism $\alpha$
from Lemma \ref{lem:homomorphisms} to $v^p+v\tilde{b'}-\tilde{b}$ yields an invariant of the
required form in $B_7$.
\end{proof}
\end{prop}

As the lfihd $\theta$ are triangular, their restriction induces lfihd on
$A_5=\kk[x,s,t,u]$, $A_6=\kk[x,y,s,t,u]$, and $A_7=\kk[x_1,x_2,x_3,y_1,y_2,y_3]$ which are
also denoted by $\theta$.
\begin{lem}\label{lem:FiniteModv}
The rings of invariants $A_5^\theta$, $A_6^\theta$,
and $A_7^\theta$ are finitely generated.
\begin{proof}
We apply the characteristic-free version (cf. \cite{rt:aacklfihd} or
\cite{hd-gk:ciagac}) of van den Essen's Algorithm (cf. \cite{ae:aacigaav}). We
do the details for (DF5), the other examples are done similarly. 
Since $\theta(s)=s+x^3U$ is a polynomial of degree $1$ in $U$ with leading
coefficient $x^3$, the invariants of the localized ring
$\kk[x,s,t,u,\frac{1}{x}]$ are generated by $1/x$, $\theta(x)|_{U=-s/x^3}=x$,
$\theta(t)|_{U=-s/x^3}=t-s^2/x^3$, and
$\theta(u)|_{U=-s/x^3}=u-3st/x^3+2s^3/x^6$.
Hence, $A_5^\theta=\kk[x,x^3t-s^2,x^6u-3x^3st+2s^3, 1/x]\cap A_5$. To obtain
generators for $A_5^\theta$, we must look at the relation ideal modulo $x$ of
the generators $f_1:=x^3t-s^2$ and $f_2:=x^6u-3x^3st+2s^3$, that is, the
preimage of the ideal $(x)$ for the map $\pi_1:\kk[X_1,X_2]\to
\kk[x,x^3t-s^2,x^6u-3x^3st+2s^3], X_1\mapsto f_1, X_2\mapsto f_2$. This is
clearly generated by $4X_1^3+X_2^2$, and so we obtain a new generator for
$\kk[x,s,t,u]^\theta$, namely
\[f_3:=\frac{1}{x^6}\pi_1(4X_1^3+X_2^2)=x^6u^2+2x^3t(2t^2-3su)+s^2(4su-3t^2).\]
Now consider $\pi_2:\kk[X_1,X_2,X_3]\to \kk[x,f_1,f_2,f_3], X_i\mapsto f_i$, $i=1,2,3$. As $\pi_2^{-1}((x))=(4X_1^3+X_2^2)\ideal \kk[X_1,X_2,X_3]$, there are no new generators. It follows that
$\kk[x,s,t,u]^\theta=\kk[x,f_1,f_2,f_3].$

For the other examples, the algorithm yields:
\begin{multline*}
A_6^\theta=\kk[x,y,x^3t-y^3s^2, x^6u-3x^3x^3st+2y^6s^3 ,\\
x^6u^2+2x^3y^3t(2t^2-3su)+y^6s^2(4su-3t^2)],\end{multline*}
and
\[A_7^\theta=\kk[x_1,x_2,x_3,x_1^3y_2-x_2^3y_1,
x_1^3y_3-x_3^3y_1,x_2^3y_3-x_3^3y_2].\qedhere\]
\end{proof}\end{lem}

We end the section with the proof of our main theorem:
\begin{proof}[Proof of Theorem \ref{thm:MainThm}]
We now show that the finite generation of the invariants modulo $v$ (Lemma
\ref{lem:FiniteModv}) together with the existence of an invariant of the form
$v^p+vb'-b$ (Proposition \ref{prop:special-invariant}) imply that the ring of
invariants is finitely generated. As the argument is the same for the three
examples, we write $B$ to denote the rings $B_5$, $B_6$, and $B_7$, and $A$ to
denote the rings $A_5$, $A_6$, and $A_7$.

As $v^p+vb'-b$ is monic as a polynomial in $v$, for any $f\in B^\theta$, there
exist unique $q$ and $r$ such that $f=q\cdot(v^p+vb'-b) +r$ and $\deg_v(r)<p$.
As $v^p+vb'-b$ and $f$ are invariant, the uniqueness of $q$ and $r$ implies
$q,r\in B^\theta$. Hence, $B^\theta$ is generated by $v^p+vb'-b$
and invariants of degree less than $p$ as polynomials in $v$.

For each degree $m<p$ the set 
$$I_m=\{ a\in A^\theta \mid a \text{ is the
leading coefficient of some } f\in B^\theta,~\deg_v(f)=m\}$$ 
is an ideal in $A^\theta$, and hence finitely generated, since $A^\theta$ is
Noetherian by Lemma \ref{lem:FiniteModv}.

Therefore, $B^\theta$ is generated by generators of $A^\theta$, $v^p+vb'-b$, and
a (finite) set of polynomials whose leading coefficients generate the ideal
$I_m$ for each $0<m<p$.
\end{proof}



\section{The 5-dimensional Example}\label{sec:dim5}

The purpose of this section is to construct, for the 5-dimensional example of
Daigle and Freudenburg (DF5), and in each positive characteristic $p$, a $w_5$-homogeneous invariant of
the form $v^p+vb'-b$, where $b,b'\in \kk[x,s,t,u]\subset B_5$ and $w_4(b)=p$, $w_4(b')=p-1$.

The proof of Theorem \ref{thm:special-invariant-B5} requires a
sequence $c_n\in \QQ[s,t,u]$, and so we must do some work over the field $\QQ$.
We consider Example (DF5) over $\QQ$ by taking $B:=\QQ[x,s,t,u,v]$ and letting $A:=\QQ[x,s,t,u]$ be
the subalgebra of $B$ with
the induced lfihd. In turn, the lfihd $\theta$ on $A$ induces a lfihd $\bar{\theta}$ on
the quotient $C:=A/(x-1)\isom \QQ[s,t,u]$ via
$\bar{\theta}(\bar{f}):=\overline{\theta(f)}$. The ring of invariants is
$C^{\bar{\theta}}=\QQ[t_1,u_1]$, where $ t_1=t-s^2$, and $u_1=u-3st+2s^3$. The grading
$w_4$ on $B$ induces a grading on $C$ which is also denoted by $w_4$.

\begin{prop}\label{prop:sequence}
Let $e:\NN\to \NN$ be given by $e(n):=\left\lfloor \frac{2n}{3}\right\rfloor$.
There exist sequences of $w_4$-homogeneous polynomials $(h_n)_{n\in\NN}$ in
$C^{\bar{\theta}}$ and $(c_n)_{n\in\NN}$ in $C$
such that $h_0=1$, $h_1=0$, 
$$c_n:=\sum_{i=0}^n \binom{n}{i} h_{n-i} s^i\in \QQ[s,t,u] \textrm{ for all }n\geq 0,$$
and for all $n\geq2$, the element $c_n$ has degree $\deg(c_n)\leq e(n)$  with respect to the standard grading $\deg$ on
$\QQ[s,t,u]$. Furthermore, for all primes $p$, the coefficients of $h_0,
h_1,\dots, h_{p-2}$ and $h_p$ are in the local ring $\ZZ_{(p)}$, and the
coefficients of $h_{p-1}$ are in $\frac{1}{p}\ZZ_{(p)}$.

\begin{rem}\label{rem:theta-cn}
For all $k,n\in\NN$,
the sequence $(c_n)_{n\in\NN}$ satisfies
$\bar{\theta}^{(k)}(c_n)=\binom{n}{k}c_{n-k}$. Indeed, we have
\begin{eqnarray*}
\bar{\theta}^{(k)}(c_n)&=& \sum_{i=0}^n \binom{n}{i} h_{n-i} \binom{i}{k}
s^{i-k} = \sum_{j=0}^{n-k} \binom{n}{j+k}h_{n-k-j} \binom{j+k}{k} s^j\\
&=& \sum_{j=0}^{n-k} \binom{n}{k}\binom{n-k}{j}h_{n-k-j}s^j=\binom{n}{k}c_{n-k}.
\end{eqnarray*}
\end{rem}

\begin{proof}[Proof of Proposition \ref{prop:sequence}]
Let $C=\bigoplus_{n\geq 0} C_n$ be the decomposition of $C$ into homogeneous
parts with respect to the grading $w_4$. For $n\geq k$, we have
$\bar{\theta}^{(k)}(C_n)\subseteq C_{n-k}$.

We will show by induction on $n$ that we can construct such sequences $h_n$ and $c_n$ so that $h_n,c_n\in C_n$. To satisfy the conditions on the coefficients of the $h_n$'s, it suffices to ensure that the denominators of the coefficients of each $h_n$ are only divisible by primes less than $n$, unless $n+1$ is a prime congruent to 1 modulo 6, in which case $n+1$ may also occur. When $n\congr 2,3,4,5 \mod{6}$, we explain how to obtain $h_n$ 
from  $\{h_j\mid 0\leq j\leq n-1\}$. The cases $n\congr 0 \mod{6}$ and $n\congr
1 \mod{6}$ must be constructed in one step, that is, when $n\congr 1
\mod{6}$, we show that we can construct $h_{n-1} $ and $h_n$ from
$\{h_j\mid0\leq j\leq n-2\}$.

Assume that $n\congr 2,3,4,5 \mod{6}$, and that
we already have $\{h_j\in C_j\mid0\leq j \leq n-1\}$ such that the
denominators of the coefficients of each $h_j$ are only divisible by primes smaller
 than $n$. By the induction
hypothesis, $c_{n-1}$ has standard degree at most $e(n-1)$ for $n>2$, and at most
$e(2)=1$
for $n=2$, that is, at most $e(n)$ in all cases. Additionally, $c_{n-1}$ is
$w_4$-homogeneous of $w_4$-degree $n-1$.
The same computation as in Remark \ref{rem:theta-cn} shows that $c:=\sum_{i=1}^n
\binom{n}{i} h_{n-i} s^i\in
\QQ[s,t,u]$ satisfies $\bar{\theta}^{(1)}(c)=nc_{n-1}$. Thus, it
will suffice  to find $h_n\in \QQ[s,t,u]^{\bar{\theta}}$ such that
$\deg(c+h_n)\leq e(n)$.

Taking 
$$c=\sum_{i+2j+3k=n} \alpha_{i,j,k} s^i t^j u^k,$$
where $\alpha_{i,j,k}\in\QQ$, one gets
\begin{multline*}
\bar{\theta}^{(1)}(c)=\sum_{i+2j+3k=n} \alpha_{i,j,k} \left( is^{i-1} t^j u^k
+ 2j  s^{i+1} t^{j-1} u^k + 3k s^i t^{j+1} u^{k-1}
\right)\\
= \sum_{i+2j+3k=n} \bigl(
(i+2)\alpha_{i+2,j-1,k} +
2j \alpha_{i,j,k} +
3(k+1)\alpha_{i+1,j-2,k+1} \bigr) s^{i+1} t^{j-1} u^k,
\end{multline*}
where $\alpha_{i',j',k'}:=0$ if $i'<0$ or $j'<0$. Since $\deg\left(
\bar{\theta}^{(1)}(c)\right)\leq e(n)$, we have 
$$(i+2)\alpha_{i+2,j-1,k} + 2j \alpha_{i,j,k} +
3(k+1)\alpha_{i+1,j-2,k+1} =0,$$
whenever $i+2j+3k=n$ and
$i+j+k\geq e(n)+1$.
Hence, each $\alpha_{i,j,k}$  is a linear combination of certain
$\alpha_{i',j',k'}$'s such that $j'<j$ and $i'+j'+k'\geq i+j+k$. Thus, if $i+j+k\geq e(n)+1$, then
$\alpha_{i,j,k}$ is a linear combination of the
$\alpha_{i',j',k'}$'s such that $j'=0$ and $i'+k'\geq e(n)+1$. Therefore, it suffices to
prove that there exists $h\in C^{\bar{\theta}}$, homogeneous of $w_4$-degree
 $n$, such that, for
all $i,k$ such that $i+3k=n$ and $i+k\geq e(n)+1$, the coefficient of $s^i u^k$ in
$h$ is equal to the coefficient of $s^i u^k$ in $c$. We then take $h_n=-h$
(note that $h$ trivially
satisfies $\deg\left(\bar{\theta}^{(1)}(h)\right)\leq e(n)$, since it is invariant, and so its coefficients  satisfy the same linear relations).

The ``relevant'' $s^i u^k$ are those where $i+3k=n$
and $i+k\geq e(n)+1$, that is, where 
$0\leq k\leq \left\lfloor \frac{n-e(n)-1}{2}\right\rfloor =
\left\lfloor\frac{n-1}{6}\right\rfloor.$ For short, we write
$d:=\lfloor\frac{n-1}{6}\rfloor$ to denote the maximum value $k$ of a relevant $s^i
u^k$.

By our assumption, $h$ is a linear combination of ``monomials'' of the form
$(-t_1)^l u_1^m$ such that $2l+3m=n$, that is, such that 
$m\in \left[0, \lfloor \frac{n}{3}\rfloor \right]$, $m\congr n \mod{2}$, and
$l=\frac{n-3m}{2}$.
Thus, for $n$ odd, we have $m=2m'+1$, where $0\leq m'\leq
\lfloor\frac{n-3}{6}\rfloor$, and for $n$ even, we have $m=2m'$ where $0\leq
m'\leq \lfloor\frac{n}{6}\rfloor$.
The coefficient of $s^i u^k$ in $(-t_1)^l u_1^m=(s^2-t)^l(u-3st+2s^3)^m$ is 
$\binom{m}{k}2^{m-k}$.

Thus, if $n\congr 2,3,4$ or $5 \mod{6}$, the number of admissible
``monomials'' in $h$ is $d+1$. Therefore, the coefficient vector $(x_0,\dots ,x_d)^T$ of $h$ is the solution of the system of linear equations
$Mx=\alpha$, where
$\alpha=\left(\alpha_{n,0,0}\, ,\dots, \alpha_{n-3d,0,d}\right)^T$, and  $M=(a_{k,m'})_{0\leq k,m'\leq d}$ is
the matrix given by
$a_{k,m'}=\binom{2m'+1}{k}2^{2m'+1-k}$, if $n$ is odd, and by
$a_{k,m'}=\binom{2m'}{k}2^{2m'-k}$, if $n$ is even.

By Lemma \ref{lem:determinants}, $M$ is invertible over
$\ZZ[\frac{1}{2}]$ in both cases. Thus, the system of equations has a unique
solution, and the primes which might occur in the denominators of the
coefficients of $h$ are $2$ and the primes occuring in the denominators of the
coefficients of some $h_j$ ($0\leq j\leq n-1$). For $n=2$, an explicit computation gives $h_2=t_1$.

If $n\congr 1\mod{6}$, the number of coefficients of $h$ is less than $d$, which
might lead to an unsolvable linear equation. Hence, we cannot use the exact same
argument. We construct $h_{n-1}$ and $h_n$ in one step. Again by induction, we
may assume that the denominators of the coefficients of $h_0,\dots, h_{n-2}$
are only divisible by primes less than $n-1$.

We want $h_{n-1}=\sum_{m'=0}^d x_{m'}(-t_1)^l u_1^{2m'}$ and 
$h_n=\sum_{m'=0}^{d-1} y_{m'}(-t_1)^l u_1^{2m'+1}$ 
such that adding $h_{n-1}$ to $\sum_{i=2}^n \binom{n-1}{i-1} h_{n-i}s^{i-1}$
cancels out the coefficient vector $\alpha$ of $\{s^{n-1}, s^{n-4}u,\dots,
s^{n-3d+2}u^{d-1} \}$, and such that adding $h_n+n h_{n-1}s$ to $\sum_{i=2}^n
\binom{n}{i} h_{n-i}s^i$ cancels out the coefficient vector $\beta$ of
$\{s^{n}, s^{n-3}u,\dots, s^{n-3d}u^{d} \}$.
Thus, we must solve $Mx=-\alpha$ and $Ny+nLx=-\beta$,
where

\begin{align*}
&M\in \Mat(d \times (d+1)), & M_{k,m'}&=\binom{2m'}{k}2^{2m'-k}, \\
&N\in \Mat((d+1) \times d), & N_{k,m'}&=\binom{2m'+1}{k}2^{2m'+1-k},\\
&L\in \Mat((d+1) \times (d+1)), & L_{k,m'}&=\binom{2m'}{k}2^{2m'-k}.
\end{align*}

Since $L$ is invertible (Lemma \ref{lem:determinants}), the two equations are
equivalent to 
$$x = -\frac{1}{n}L^{-1}(\beta+Ny) \quad \text{ and }\quad
ML^{-1}Ny= -ML^{-1}\beta+n \alpha.$$
As $ML^{-1}N$ is the $d\times d$ matrix with entries
$\binom{2m'+1}{k}2^{2m'+1-k}$, by Lemma \ref{lem:determinants}, it is invertible.

Therefore, we can find unique vectors $x$ and $y$. From the equations giving $x$
and $y$, we deduce that
the denominators of the entries of $y$ are only divisible by primes less than
$n-1$, and if $n$ is prime, the denominators of the entries of $x$ may have the additional prime factor $n$.
\end{proof}\end{prop}

\begin{lem}\label{lem:determinants}
 Let $d\in\NN$. If $M_e, M_o\in\Mat((d+1)\times (d+1),\ZZ)$ are given by
$$(M_e)_{k,m}:= \binom{2m}{k}2^{2m-k}\quad \text{and} \quad 
(M_o)_{k,m}:= \binom{2m+1}{k}2^{2m+1-k}$$
for all $k,m=0,\dots, d$, then
$ \det(M_e)=2^{d(d+1)}$, 
and $\det(M_o)=2^{(d+1)^2}$.
\begin{proof}
 If $x_0,\dots, x_d$ are indeterminates, denote by $V(x_0,\dots, x_d)$ the Vandermonde
matrix
$$V(x_0,\dots, x_d)= \begin{pmatrix}  1 & x_0 & \dots  &   x_0^d  \\
				1 & x_1 & \dots  &   x_1^d  \\
				\vdots & \vdots & \ddots & \vdots  \\
				1 & x_d & \dots  &   x_d^d \end{pmatrix}.$$
By definition of the binomial coefficients, we have
$$V(x_0,\dots, x_d)\cdot M_e \congr V((2+x_0)^2 ,\dots, (2+x_d)^2)
\mod{(x_0^{d+1},\dots, x_d^{d+1})},$$
and therefore  their determinants are also congruent modulo $(x_0^{d+1},\dots,
x_d^{d+1})$.
Since  $\det\left(V(x_0,\dots, x_d)\right)=\prod_{j>i} (x_j-x_i)$, and
\begin{eqnarray*}\det\left(V((2+x_0)^2 ,\dots, (2+x_d)^2)\right)&=&\prod_{j>i}
\left((2+x_j)^2-(2+x_i)^2\right)\\
&=& \prod_{j>i} (x_j-x_i) \prod_{j>i} (4+x_j+x_i),
\end{eqnarray*} 
we obtain
$$\prod_{j>i} (x_j-x_i)\cdot \left[ \det(M_e)- \prod_{j>i}
(4+x_j+x_i)\right]\congr 0 \mod{(x_0^{d+1},\dots, x_d^{d+1})}.$$
This implies that the coefficient of $x_1x_2^2\cdots x_d^d$ of the left hand
side, namely $\det(M_e)- \prod_{j>i} 4$, has to vanish. Hence,
$\det(M_e)=2^{d(d+1)}$.

A similar argument proves the statement about $M_o$. The key is to recognize that
$V(x_0,\dots, x_d)\cdot M_o$ is congruent 
to 
$$\begin{pmatrix}
 2+x_0 & (2+x_0)^3 & \dots & (2+x_0)^{2d+1} \\
 2+x_1 & (2+x_1)^3 & \dots & (2+x_1)^{2d+1} \\
\vdots & \vdots & \ddots & \vdots  \\
2+x_d & (2+x_d)^3 & \dots & (2+x_d)^{2d+1}
\end{pmatrix}$$
modulo $(x_0^{d+1},\dots, x_d^{d+1})$, and that this matrix has determinant equal to
$\prod_{l=0}^d (2+x_l)\cdot \det\left(V((2+x_0)^2 ,\dots, (2+x_d)^2)\right)$.
\end{proof}\end{lem}

We are now prepared to prove the existence of the special invariant.

\begin{thm}\label{thm:special-invariant-B5}
Let $\kk$ be of positive characteristic $p$, and let $(B_5,\theta)$ be
as in Example (DF5) over $\kk$.\\ There exists a $w_5$-homogeneous invariant in
$B_5$ of the form $v^p+vb'-b$, where $b',b\in \kk[x,s,t,u]\subset
B_5$ have $w_4$-degree
$w_4(b)=p$ and $w_4(b')=p-1$.
\begin{proof}
If suffices to find $b,b'\in \kk[x,s,t,u]$ such that
$\theta(b)=b+x^2b'U+x^{2p}U^p$. Indeed, this implies $\theta(b')=b'$, and
\begin{multline*}
\theta(v^p+vb'-b)=\theta(v)^p+\theta(v)\theta(b')-\theta(b)\\
=(v+x^2U)^p+(v+x^2U )b'-(b+x^2b'U+x^{2p}U^p) = v^p+vb'-b.
\end{multline*}
 Let $\O:=\WW(\kk)$ be the Witt ring of $\kk$, and let $\KK$ be
the field of fractions of $\O$. Hence $\KK$ is a discrete valued field of
unequal
characteristic with valuation ring $\O$, valuation ideal $(p)$ and residue
class field $\kk$.  Example (DF5) over $\KK$ has a lfihd which restricts to
$\O[x,s,t,u,v]$. Reduction modulo $p$ then leads to Example (DF5) over
$\kk$. Thus, to obtain $b,b'\in \kk[x,s,t,u]$ such that
$\theta(b)=b+x^2b'U+x^{2p}U^p$, it
suffices to find $\tilde{b},\tilde{b}'\in \O[x,s,t,u]\subseteq \KK[x,s,t,u]$
such that $\theta(\tilde{b})
\congr \tilde{b}+x^2\tilde{b}'U+x^{2p}U^p \mod{p}.$

Let $(c_n)_{n\in\NN}\subset \QQ[s,t,u]\subseteq \KK[s,t,u]$ be
the sequence constructed in Proposition \ref{prop:sequence}. Set
$$b_n:=x^{2p} c_n(\frac{s}{x^3},\frac{t}{x^3},\frac{u}{x^3}) \text{ for } 0\leq
n\leq p,$$
that is, $b_n$ is the homogenization of $c_n$ of degree $2p$ with respect
to the
grading $w_5$. By construction of
$c_n$, the elements $b_n$ are indeed in $\KK[x,s,t,u]$ and have coefficients
in $\ZZ_{(p)}\subseteq \O$ resp. $\frac{1}{p}\ZZ_{(p)}$ for $n=p-1$. Moreover,
we have $b_0=x^{2p}$, and if $p> 2$, $x^2$ divides $b_{p-1}$. A
similar calculation as in Remark \ref{rem:theta-cn} shows that $\theta^{(k)}(b_p)=\binom{p}{k}b_{p-k}$ for all $0\leq k\leq p$, and $\theta^{(k)}(b_p)=0$ for $k>p$.
Hence, if $p=2$, $\theta(b_2)=b_2+x^4U^2 \mod{2}$, and otherwise, $\theta(b_p)\congr b_p+x^2\left(\frac{pb_{p-1}}{x^2}\right)U+x^{2p}U^p
\mod{p}$, as desired.
\end{proof}
\end{thm}


\section{Generalized Form of the Examples}\label{sec:generalization}

In this section we explain how the arguments of Sections \ref{sec:main-results}
and \ref{sec:dim5} can be adapted to generalized forms of examples (R7), (F6),
and (DF5). We start by writing down these general forms. In characteristic
zero, these generalizations
are  discussed in Section~7.2.3 of Freudenburg's
book (c.f.~\cite{gf:atlnd}). Note that the version of (R7) presented here is the
original version. Let $m\geq 2$ be an integer.

\noindent 
{\bf Generalized Daigle-Freudenburg's example} (DF5-m): \\$B_{5}:=\kk[x,s,t,u,v]$,
\begin{align*}
\theta(x)&=x,& \theta(s)&=s+x^{m+1} U, \\
\theta(t)&=t+2sU+x^{m+1}U^2,& \theta(u)&=u+3tU+3sU^2+x^{m+1}U^3,\\ 
\theta(v)&=v+x^m U. && 
\end{align*} 
{\bf Generalized Freudenburg's example} (F6-m): \\$B_{6}:=\kk[x,y,s,t,u,v]$,
\begin{align*} 
\theta(x)&=x,&    \theta(t)&=t+2y^{m+1}sU+x^{m+1}y^{m+1}U^2,\\
\theta(y)&=y,& \theta(u)&=u+3y^{m+1}tU+3y^{2(m+1)}sU^2+x^{m+1}y^{2(m+1)}U^3,\\
\theta(s)&=s+x^{m+1} U,& \theta(v)&=v+x^my^m U.
\end{align*}
{\bf Generalized Roberts's example} (R7-m):\\ $B_{7}:=\kk[x_1,x_2,x_3,y_1,y_2,y_3,v]$,
\begin{align*}
\theta(x_i)&=x_i,& \theta(y_i)&=y_i+x_i^{m+1} U\quad (i=1,2,3),\\ 
\theta(v)&=v+ x_1^m x_2^m x_3^m U.&&
\end{align*}

On the $\kk$-algebras $B_{5}, B_{6}, B_{7}$, we define gradings $w_5,w_6,w_7$ by assigning the following degrees:
\[\begin{array}{l}
w_5(x)=1,~w_5(s)=w_5(t)=w_5(u)=m+1,~ w_5(v)=m;\\
w_6(x)=w_6(y)=1,~w_6(s)=m+1,~w_6(t)=2(m+1),\\
\hspace*{3cm} w_6(u)=3(m+1), w_6(v)=2m;\\
w_7(x_i)=1,~w_7(y_i)=m+1,~(i=1,2,3),~w_7(v)=3m.\\
  \end{array}\]
With respect to these gradings, the lfihd and the corresponding $\Ga$-actions
are homogeneous, and so the rings of invariants are graded subalgebras. We will also use an additional grading
$w_4$ on $B_5$ which is given by:
$$w_4(x)=0,~w_4(s)=1,~w_4(t)=2,~w_4(u)=3,~w_4(v)=1.$$

\begin{prop}\label{prop:special-invariant-m}
 In each of $B_5$, $B_6$, and $B_7$, there exists a homogeneous invariant of the
form $v^p+vb'-b$ such that $v$ does not appear in $b$ and $b'$.
\begin{proof}
A short computation verifies that there exists similar relationship between the generalized examples as establish in Lemma \ref{lem:homomorphisms}. Namely, the ring $B_{5}$ is isomorphic to $B_{6}/(y-1)$ and the lfihd on $B_{5}$ is the lfihd induced by this isomorphism from the lfihd on $B_{6}$, and a homomorphism $\alpha_m:B_6\to B_7$ respecting the lfihd is
given by:
\begin{align*}
\alpha(x)&=x_1,& \alpha(y)&=x_2 x_3,& \alpha(s)&=y_1,\\
\alpha(t)&\multispan{3}$\,=(x_3^{m+1}y_2+x_2^{m+1}y_3)y_1-x_1^{m+1}
y_2y_3,$& \alpha(v)&=v, \\
\alpha(u)&\multispan{5}$\,=(x_2^{2(m+1)}y_3^2+x_2^{m+1}x_3^{m+1}y_2y_3+x_3^{
2(m+1) }y_2^2)y_1-
(x_3^{m+1}y_2+x_2^{m+1}y_3)x_1^{m+1}y_2y_3$
\end{align*}

As in Proposition \ref{prop:special-invariant}, we construct the special
invariant for Example (DF5-m), and then the relationships described above yields
the special invariants for the two other examples.

To construct the special invariant for (DF5-m), it suffices as in Theorem
\ref{thm:special-invariant-B5} to find $b,b'\in \kk[x,s,t,u]$ such that
$\theta(b)=b+x^mb'U+x^{mp}U^p$. Similarly, it suffices to find $\tilde{b},
\tilde{b}'\in \O[x,s,t,u] \subseteq \KK[x,s,t,u]$ such that $\theta(\tilde{b})
\congr \tilde{b}+x^m \tilde{b}'U+x^{mp}U^p \mod{p}$. Let
$(c_n)_{n\in\NN}\subset \QQ[s,t,u]\subseteq \KK[s,t,u]$ be
the sequence constructed in Proposition \ref{prop:sequence}. Set
$$b_n:=x^{mp} c_n(\frac{s}{x^{m+1}},\frac{t}{x^{m+1}},\frac{u}{x^{m+1}}) \text{
for } 0\leq n\leq p,$$
that is, $b_n$ is the homogenization of $c_n$ of degree $mp$ with respect to the
grading $w_5$. By construction of $c_n$, the elements $b_n$ are indeed in
$\KK[x,s,t,u]$ and have coefficients in $\ZZ_{(p)}\subseteq \O$ resp.
$\frac{1}{p}\ZZ_{(p)}$ for $n=p-1$. Moreover, we have $b_0=x^{mp}$, and if
$p>2$, $x^m$ divides $b_{p-1}$. A similar calculation
as in Remark \ref{rem:theta-cn} shows that
$\theta^{(k)}(b_p)=\binom{p}{k}b_{p-k}$ for all $0\leq k\leq p$, and
$\theta^{(k)}(b_p)=0$ for $k>p$. 
Hence, if $p=2$, $\theta(b_2)=b_2+x^{2m}U^2 \mod{2}$, and otherwise,
$\theta(b_p)\congr b_p+x^m\left(\frac{pb_{p-1}}{x^m}\right)U+x^{mp}U^p
\mod{p}$, as desired.

\end{proof}\end{prop}

This finally leads to the generalized version of our main theorem:

\begin{thm}\label{thm:MainThm-m}
In every positive characteristic, the rings of invariants $B_5^\theta$,
$B_6^\theta$, and $B_7^\theta$, are finitely generated.
\begin{proof}
 The argument is as in the proof of Theorem \ref{thm:MainThm}. Again, the
$\kk$-algebra $A$ (that is, $A_5$, $A_6$ and $A_7$) is finitely generated, and
the finite generation of $B^\theta$ then follows from the special invariant
established by Proposition \ref{prop:special-invariant-m}.
\end{proof}
\end{thm}

\begin{rem}~\
\begin{enumerate}
\item For (R7-m), our argument has the additional advantage of providing a
somewhat simpler proof than Kurano's, admittedly more general, argument
(c.f.~\cite{kk:pcfgsrarcfph}). Moreover, our proof constructs an invariant which
is monic of degree $p$ as a polynomial in $v$, where Kurano proves only the
existence of some invariant monic as a polynomial in $v$.

\item In Section~7.2.3 of his book (c.f.~\cite{gf:atlnd}), Freudenburg explains
how (DF5-m) is really (R7-m) with all the symmetries removed. Precisely he
proves that $B_{5}^\theta$ is equal to the ring of invariants of $B_{7}^\theta$
under the action of a reductive group. In characteristic zero, this
automatically implies that $B_{5}^\theta$ is not finitely generated.
Freudenburg's argument concerning the relationship between $B_{7}^\theta$ and
$B_{5}^\theta$ remains valid in positive characteristic, and so the finite
generation of $B_{7}^\theta$ implies the finite generation of $B_{5}^\theta$. In
characteristic zero, since $B_{6}^\theta$ surjects onto $B_{5}^\theta$, it
follows that $B_{6}^\theta$ is also not finitely generated. It appears that, in
positive characteristic, the use of our constructive argument is unavoidable to
prove the finite generation of $B_{6}^\theta$. \end{enumerate}
\end{rem}


\bibliographystyle{plain}
\bibliography{reference}

\vspace*{.5cm}
\end{document}